\documentclass[11pt, twoside]{article}
\usepackage{latexsym}
\usepackage{amsmath}
\usepackage{amssymb}
\usepackage[all]{xy}
\usepackage{amsfonts}
\usepackage{verbatim}
\usepackage{amsthm}
\usepackage{bm}
\usepackage{mathrsfs}
\usepackage{epsfig}
\usepackage{xy}
\usepackage{array}
\usepackage{stmaryrd}
\usepackage{graphicx,color}
\usepackage{xcolor}
\usepackage{tikz}
\usetikzlibrary{arrows,calc}
\usepackage{etex}
\usepackage{mathdots}
\usepackage{float}
\usepackage{graphics}
\usepackage{pdflscape}
\usepackage{CJK}
\usepackage{anysize,hyperref}
\input xypic
\xyoption{all}
\usepackage{extarrows}
\usepackage[perpage,symbol]{footmisc}
\usepackage{setspace}
\def\C{\mathcal{C}}

\def\dr{\ar@{->}[r]}

\def\X{\mathscr{X}}

\def\Hom{\mbox{Hom}}

\def\rad{\mbox{rad}}
\def\End{\mbox{End}}
\def\Ker{\mathsf{Ker}\hspace{.01in}}
\def\Coker{\mathsf{Coker}\hspace{.01in}}

\newcommand{\ind}{\mathsf{ind}\hspace{.01in}}
\newcommand{\supp}{\mathsf{Supp}\hspace{.01in}}
\begin{document}
\baselineskip=15pt
\title{\Large{\bf Grothendieck groups and Auslander-Reiten $\bm{(d+2)}$-angles$^\bigstar$\footnotetext{\hspace{-1em}$^\bigstar$This work was supported by the National Natural Science Foundation of China (Grant No. 11901190 and 11671221), and by the Hunan Provincial Natural Science Foundation of China (Grant No. 2018JJ3205), and by the Scientific Research Fund of Hunan Provincial Education Department (Grant No. 19B239).}}}
\medskip
\author{Panyue Zhou}

\date{}

\maketitle
\def\blue{\color{blue}}
\def\red{\color{red}}

\newtheorem{theorem}{Theorem}[section]
\newtheorem{lemma}[theorem]{Lemma}
\newtheorem{corollary}[theorem]{Corollary}
\newtheorem{proposition}[theorem]{Proposition}
\newtheorem{conjecture}{Conjecture}
\theoremstyle{definition}
\newtheorem{definition}[theorem]{Definition}
\newtheorem{question}[theorem]{Question}
\newtheorem{remark}[theorem]{Remark}
\newtheorem{remark*}[]{Remark}
\newtheorem{example}[theorem]{Example}
\newtheorem{example*}[]{Example}
\newtheorem{condition}[theorem]{Condition}
\newtheorem{condition*}[]{Condition}
\newtheorem{construction}[theorem]{Construction}
\newtheorem{construction*}[]{Construction}

\newtheorem{assumption}[theorem]{Assumption}
\newtheorem{assumption*}[]{Assumption}

\baselineskip=17pt
\parindent=0.5cm

\begin{abstract}
\baselineskip=16pt
Xiao and Zhu has shown that if $\C$ is a locally finite triangulated category, then
the Auslander-Reiten triangles generate the relations for the Grothendieck group of $\C$.
The notion of $(d +2)$-angulated categories is a ``higher dimensional" analogue
of triangulated categories. In this article, we show that if A $(d+2)$-angulated category $\C$ is locally finite if and only if
the Auslander-Reiten $(d+2)$-angles generate the relations for the Grothendieck group of $\C$.
This extends the result of Xiao and Zhu, and gives the
converse of Xiao and Zhu's result is also true.\\[0.5cm]
\medskip
\end{abstract}

\pagestyle{myheadings}
\markboth{\rightline {\scriptsize   Panyue Zhou}}
         {\leftline{\scriptsize  Grothendieck groups and Auslander-Reiten $(d+2)$-angles}}

\section{Introduction}
Auslander-Reiten theory was introduced by Auslander and Reiten in \cite{ar1,ar2}. Since its introduction, Auslander-Reiten theory has become a fundamental tool for
studying the representation theory of Artin algebras. It is well-known that
a module category of over an  Artin algebra has almost split sequences.
For an Artin algebra of finite type, Butler \cite{bu} proved that
the relations of its Grothendieck group are generated by all Auslander-Reiten
sequences. Soon later Auslander showed that the converse is true
in \cite{au}. The notion of Auslander-Reiten triangles in a triangulated category was
introduced by Happel in \cite{ha}. In contrast to  module categories over Artin algebras, not all triangulated
categories  have
Auslander-Reiten triangles \cite{ha}. It was proved in \cite{ha} that the derived category
of a finite dimensional algebra has Auslander-Reiten triangles if and only if the
global dimension of the algebra is finite.
Reiten and Van den Bergh \cite{rv}
proved that the existence of Auslander-Reiten triangles if and only if the
existence of Serre functor in a triangulated category. Recently, Xiao and Zhu \cite{xz}
showed that if $\C$ is a locally finite triangulated category, then
the Auslander-Reiten triangles generate the relations for the Grothendieck group of $\C$. Beligiannis \cite{be} proved the converse of this result holds when $\C$ is a compactly generated triangulated category. In more recent times, many authors have shown the reverse direction of
Xiao and Zhu is true in some special cases \cite{h1,pppp}.
Extriangulated categories were introduced by Nakaoka and Palu \cite{np} as a simultaneous
generalization of exact categories and triangulated categories.
Hence, many results hold on exact categories and triangulated categories can be unified in the same framework.
Iyama, Nakaoka and Palu \cite{inp}  introduced the notion of almost split extensions and Auslander-Reiten-Serre duality for extriangulated categories, and gave explicit connections between these notions and also with the
classical notion of dualizing $k$-varieties.
Zhu and Zhuang \cite{zz} has shown that a locally finite  extriangulated category $\C$ has Auslander-Reiten $\mathbb{E}$-triangles and
the relations of Grothendieck group are generated by the Auslander-Reiten $\mathbb{E}$-triangles. A partial converse result is given when restricting to a triangulated category with a cluster tilting subcategory.

In \cite{gko}, Geiss, Keller and Oppermann introduced $(d+2)$-angulated categories. These are
are a ``higher dimensional" analogue of triangulated categories,
 in the sense that triangles are replaced by $(d+2)$-angles, that is, morphism sequences of length $(d+2)$. Thus a $1$-angulated category is precisely a triangulated category.
 They appear for example as certain cluster tilting subcategories of triangulated categories.
Iyama and Yoshino defined the notion of Auslander-Reiten $(d+2)$-angles in special $(d+2)$-angulated categories.
This notion was generalized to arbitrary $(d+2)$-angulated categories by Fedele \cite{fe2}.
Fedele also proved that there are Auslander-Reiten $(d+2)$-angles in certain subcategories of $(d+2)$-angulated categories.
The author \cite{z1} showed that if a $(d+2)$-angulated category has Auslander-Reiten $(d+2)$-angles if and only if $\C$ has a Serre functor.
Moreover, the author \cite{z2} also proved if $\C$ is a locally $(d+2)$-angulated category, then $\C$ has Auslander-Reiten $(d+2)$-angles.

Bergh and Thaule \cite{bt2} defined the Grothendieck group of a $(d+2)$-angulated category. As in the triangulated case, it
is the free abelian group on the set of isomorphism classes of objects, modulo the Euler relations corresponding to the
$(d+2)$-angles. They showed that when $d$ is odd, the set of subgroups corresponds bijectively to the complete and dense
$(d+2)$-angulated subcategories. Fedele \cite{fe1} showed that under suitable circumstances, the Grothendieck group of a triangulated category can be expressed as a quotient of the split Grothendieck group of a higher cluster tilting subcategory of a triangulated category.
Recently, Herschend, Liu and Nakaoka \cite{hln} defined $d$-exangulated categories as a ``higher dimensional" analogue
of extriangulated categories. Many categories studied in representation
theory turn out to be $d$-exangulated. In particular, $d$-exangulated categories
simultaneously generalize $(d+2)$-angulated and $d$-exact categories \cite{ja}.
Haugland \cite{h2} defined the Grothendieck group of a $d$-exangulated category, and
classified dense complete subcategories of $d$-exangulated categories with a $d$-(co)generator in terms of subgroups of the Grothendieck group.

The aim of this article to discuss the relation between Grothendieck groups in $(d+2)$-angulated categories and
Auslander-Reiten $(d+2)$-angles in $(d+2)$-angulated categories. We show that A $(d+2)$-angulated category is locally finite if and only if
the Auslander-Reiten $(d+2)$-angles generate the relations for the Grothendieck group of this $(d+2)$-angulated category, see Theorem \ref{main1} and Theorem \ref{main2}.
This extends the result of Xiao and Zhu, and gives the
converse of Xiao and Zhu's result holds.
We hope that our work would motivate and
further study on $(d+2)$-angulated categories.

This article is organised as follows: In Section 2, we review some elementary definitions
that we need to use, including  $(d+2)$-angulated categories and Auslander-Reiten $(d+2)$ angles.
In Section 3, we show our two main results.

\section{Preliminaries}
In this section, we first recall the definition and basic properties of $(d+2)$-angulated categories from \cite{gko}.
Let $\C$ be an additive category with an automorphism $\Sigma^d:\C\rightarrow\C$, and an integer $d$ greater than or equal to one.

A $(d+2)$-$\Sigma^d$-$sequence$ in $\C$ is a sequence of objects and morphisms
$$A_0\xrightarrow{f_0}A_1\xrightarrow{f_1}A_2\xrightarrow{f_2}\cdots\xrightarrow{f_{d-1}}A_n\xrightarrow{f_d}A_{d+1}\xrightarrow{f_{d+1}}\Sigma^d A_0.$$
Its {\em left rotation} is the $(d+2)$-$\Sigma^d$-sequence
$$A_1\xrightarrow{f_1}A_2\xrightarrow{f_2}A_3\xrightarrow{f_3}\cdots\xrightarrow{f_{d}}A_{d+1}\xrightarrow{f_{d+1}}\Sigma^d A_0\xrightarrow{(-1)^{d}\Sigma^d f_0}\Sigma^d A_1.$$
A \emph{morphism} of $(d+2)$-$\Sigma^d$-sequences is  a sequence of morphisms $\varphi=(\varphi_0,\varphi_1,\cdots,\varphi_{d+1})$ such that the following diagram commutes
$$\xymatrix{
A_0 \ar[r]^{f_0}\ar[d]^{\varphi_0} & A_1 \ar[r]^{f_1}\ar[d]^{\varphi_1} & A_2 \ar[r]^{f_2}\ar[d]^{\varphi_2} & \cdots \ar[r]^{f_{d}}& A_{d+1} \ar[r]^{f_{d+1}}\ar[d]^{\varphi_{d+1}} & \Sigma^d A_0 \ar[d]^{\Sigma^d \varphi_0}\\
B_0 \ar[r]^{g_0} & B_1 \ar[r]^{g_1} & B_2 \ar[r]^{g_2} & \cdots \ar[r]^{g_{d}}& B_{d+1} \ar[r]^{g_{d+1}}& \Sigma^d B_0
}$$
where each row is a $(d+2)$-$\Sigma^d$-sequence. It is an {\em isomorphism} if $\varphi_0, \varphi_1, \varphi_2, \cdots, \varphi_{d+1}$ are all isomorphisms in $\C$.

\begin{definition}\cite[Definition 2.1]{gko}
A $(d+2)$-\emph{angulated category} is a triple $(\C, \Sigma^d, \Theta)$, where $\C$ is an additive category, $\Sigma^d$ is an automorphism of $\C$ ($\Sigma^d$ is called the $d$-suspension functor), and $\Theta$ is a class of $(d+2)$-$\Sigma^d$-sequences (whose elements are called $(d+2)$-angles), which satisfies the following axioms:
\begin{itemize}
\item[\textbf{(N1)}]
\begin{itemize}
\item[(a)] The class $\Theta$ is closed under isomorphisms, direct sums and direct summands.

\item[(b)] For each object $A\in\C$ the trivial sequence
$$ A\xrightarrow{1_A}A\rightarrow 0\rightarrow0\rightarrow\cdots\rightarrow 0\rightarrow \Sigma^dA$$
belongs to $\Theta$.

\item[(c)] Each morphism $f_0\colon A_0\rightarrow A_1$ in $\C$ can be extended to $(d+2)$-$\Sigma^d$-sequence: $$A_0\xrightarrow{f_0}A_1\xrightarrow{f_1}A_2\xrightarrow{f_2}\cdots\xrightarrow{f_{d-1}}A_d\xrightarrow{f_d}A_{d+1}\xrightarrow{f_{d+1}}\Sigma^d A_0.$$
\end{itemize}
\item[\textbf{(N2)}] A $(d+2)$-$\Sigma^n$-sequence belongs to $\Theta$ if and only if its left rotation belongs to $\Theta$.

\item[\textbf{(N3)}] Each solid commutative diagram
$$\xymatrix{
A_0 \ar[r]^{f_0}\ar[d]^{\varphi_0} & A_1 \ar[r]^{f_1}\ar[d]^{\varphi_1} & A_2 \ar[r]^{f_2}\ar@{-->}[d]^{\varphi_2} & \cdots \ar[r]^{f_{d}}& A_{d+1} \ar[r]^{f_{d+1}}\ar@{-->}[d]^{\varphi_{d+1}} & \Sigma^d A_0 \ar[d]^{\Sigma^d \varphi_0}\\
B_0 \ar[r]^{g_0} & B_1 \ar[r]^{g_1} & B_2 \ar[r]^{g_2} & \cdots \ar[r]^{g_{d}}& B_{d+1} \ar[r]^{g_{d+1}}& \Sigma^d B_0
}$$ with rows in $\Theta$, the dotted morphisms exist and give a morphism of  $(d+2)$-$\Sigma^d$-sequences.

\item[\textbf{(N4)}] In the situation of (N3), the morphisms $\varphi_2,\varphi_3,\cdots,\varphi_{d+1}$ can be chosen such that the mapping cone
$$A_1\oplus B_0\xrightarrow{\left(\begin{smallmatrix}
                                        -f_1&0\\
                                        \varphi_1&g_0
                                       \end{smallmatrix}
                                     \right)}
A_2\oplus B_1\xrightarrow{\left(\begin{smallmatrix}
                                        -f_2&0\\
                                        \varphi_2&g_1
                                       \end{smallmatrix}
                                     \right)}\cdots\xrightarrow{\left(\begin{smallmatrix}
                                        -f_{d+1}&0\\
                                        \varphi_{d+1}&g_d
                                       \end{smallmatrix}
                                     \right)} \Sigma^n A_0\oplus B_{d+1}\xrightarrow{\left(\begin{smallmatrix}
                                        -\Sigma^d f_0&0\\
                                        \Sigma^d\varphi_1&g_{d+1}
                                       \end{smallmatrix}
                                     \right)}\Sigma^dA_1\oplus\Sigma^d B_0$$
belongs to $\Theta$.
   \end{itemize}
\end{definition}

The following two lemmas are very useful which are needed later on.

\begin{lemma}\emph{\cite[Lemma 3.13]{fe2}}\label{y1}
Let $\C$ be a $(d+2)$-angulated category, and
\begin{equation}\label{t1}
\begin{array}{l}
A_0\xrightarrow{\alpha_0}A_1\xrightarrow{\alpha_1}A_2\xrightarrow{\alpha_2}\cdots\xrightarrow{\alpha_{d-1}}A_d\xrightarrow{\alpha_d}A_{d+1}\xrightarrow{\alpha_{d+1}}\Sigma^d A_0.\end{array}
\end{equation}
a $(d+2)$-angle in $\C$. Then the following are equivalent:
\begin{itemize}
\item[\rm (1)] $\alpha_0$ is a section (also known as a split monomorphism);

\item[\rm (2)] $\alpha_d$ is a retraction (also known as a split epimorphism);

\item[\rm (3)] $\alpha_{d+1}=0$.
\end{itemize}
If a $(d+2)$-angle \emph{(\ref{t1})} satisfies one of the above equivalent conditions, it is called \emph{split}.
\end{lemma}

\begin{lemma}\emph{\label{m4}\cite[Lemma 4.1]{bt1}}
Suppose $\C$ is a $(d+2)$-angulated category, and let
$$\xymatrix{
A_0\ar[r]^{\alpha_1}& A_1 \ar[r]^{\alpha_2} & A_2 \ar[r]^{\alpha_3}& \cdots \ar[r]&A_{d-1}\ar[r]^{\alpha_{d-1}}&A_{d}\ar[r]^{\alpha_d}\ar[d]^{\varphi_d}&A_{d+1}\ar[r]^{\alpha_{d+1}}\ar@{=}[d] & \Sigma^d A_0 \\
B_0\ar[r]^{\beta_0}&B_1 \ar[r]^{\beta_2} & B_3 \ar[r]^{\beta_3}& \cdots \ar[r]&B_{d-1}\ar[r]^{\beta_{d-1}} & B_{d} \ar[r]^{\beta_d}&A_{d+1}\ar[r]^{\beta_{d+1}}& \Sigma^dB_0}$$
be commutative diagram whose rows are $(d+2)$-angles.
Then it can be completed the diagram to a morphism
$$\xymatrix{
A_0\ar[r]^{\alpha_0}\ar@{-->}[d]^{\varphi_{0}}& A_1 \ar[r]^{\alpha_1}\ar@{-->}[d]^{\varphi_{1}} & A_2 \ar[r]^{\alpha_2}\ar@{-->}[d]^{\varphi_{2}}& \cdots \ar[r]&A_{d-1}\ar[r]^{\alpha_{d-1}}\ar@{-->}[d]^{\varphi_{d-1}}&A_{d}\ar[r]^{\alpha_d}\ar[d]^{\varphi_d}&A_{d+1}\ar[r]^{\alpha_{d+1}}\ar@{=}[d] & \Sigma^d A_0 \ar@{-->}[d]^{\Sigma^d\varphi_0}\\
B_0\ar[r]^{\beta_0}&B_1 \ar[r]^{\beta_1} & B_3 \ar[r]^{\beta_2}& \cdots \ar[r]&B_{d-1}\ar[r]^{\beta_{d-1}} & B_{d} \ar[r]^{\beta_d}&A_{d+1}\ar[r]^{\beta_{d+1}}& \Sigma^dB_0}$$
of $(d+2)$-angles such that
$$A_0\xrightarrow{\left(
                    \begin{smallmatrix}
                      -\alpha_1 \\
                      \varphi_0 \\
                    \end{smallmatrix}
                  \right)} A_1\oplus B_0\xrightarrow{\left(
                             \begin{smallmatrix}
                               \alpha_1 & 0 \\
                               \varphi_1 & \beta_0
                             \end{smallmatrix}
                           \right)}
 A_2\oplus B_1\xrightarrow{\left(
                                       \begin{smallmatrix}
                                         -\alpha_2 & 0  \\
                                         \varphi_2 &\beta_1
                                       \end{smallmatrix}
                                     \right)}\cdots\xrightarrow{\left(
                                       \begin{smallmatrix}
                                         (-1)^{d}\alpha_{d-1}& 0  \\
                                         \varphi_{d-1} &\beta_{d-2}
                                       \end{smallmatrix}
                                     \right)}A_d\oplus B_{d-1}$$
                                     $\xrightarrow{\left(
                                       \begin{smallmatrix}
                                         (-1)^{d+1}\varphi_d,&\beta_{d-1}
                                       \end{smallmatrix}
                                     \right)}B_d\xrightarrow{~\alpha_{d+1}\beta_d~}\Sigma^dA_0$
 is a $(d+2)$-angle in $\C$.
\end{lemma}

\medskip
Now we recall an Auslander-Reiten $(d+2)$-angle in $(d+2)$-angulated categories.

We denote by ${\rm rad}_{\C}$ the Jacobson radical of $\C$. Namely, ${\rm rad}_{\C}$ is an ideal of $\C$ such that ${\rm rad}_{\C}(A, A)$
coincides with the Jacobson radical of the endomorphism ring ${\rm End}(A)$ for
any $A\in\C$.

\begin{definition}\cite[Definition 3.8]{iy} and \cite[Definition 5.1]{fe2}
Let $\C$ be a $(d+2)$-angulated category.
A $(d+2)$-angle
$$A_{\bullet}:~~\xymatrix {A_0 \xrightarrow{~\alpha_0~}A_1 \xrightarrow{~\alpha_1~} A_2 \xrightarrow{~\alpha_2~} \cdots
  \xrightarrow{~\alpha_{d - 1}~} A_d \xrightarrow{~\alpha_{d}~} A_{d+1}\xrightarrow{~\alpha_{d+1}~} \Sigma^d A_0}$$
in $\C$ is called an \emph{Auslander-Reiten $(d+2)$-angle }if
$\alpha_0$ is left almost split, $\alpha_d$ is right almost split and
when $d\geq 2$, also $\alpha_1,\alpha_2,\cdots,\alpha_{d-1}$ are in $\rad_{\C}$.
\end{definition}

\begin{remark}\cite[Remark 5.2]{fe2}
Assume $A_{\bullet}$ as in the above definition is an Auslander-Reiten $(d+2)$-angle.
Since $\alpha_0$ is left almost split implies that $\End(A_0)$ is local and hence $A_0$ is
indecomposable. Similarly, since $\alpha_d$ is right almost split, then $\End(A_{d+1})$ is local and hence $A_{d+1}$ is indecomposable.
Moreover, when $d=1$, we have $\alpha_0$ and $\alpha_d$ in $\rad_{\C}$, so that $\alpha_d$ is right minimal and $\alpha_0$ is left minimal. When $d\geq 2$, since $\alpha_{d-1}\in\rad_{\C}$, we have that $\alpha_d$ is right minimal and similarly $\alpha_0$ is left
minimal.
\end{remark}

Let $k$ be an algebraically closed
field and $\C$ be a $k$-linear Hom-finite additive category.  Recall that the notion of a Serre functor.
We call a $k$-linear autofunctor $\mathbb{S}\colon \C\to\C$ a \emph{Serre functor} of $\C$  if there exists a functorial isomorphism:
$$\Hom_{\C}(X,Y)\simeq D\Hom_{\C}(Y, \mathbb{S}X)$$
for any $X,Y\in\C$, where $D(-)=\Hom_k(-, k)$ is the $k$-linear duality functor.

The author gave a relationship between an Auslander-Reiten $(d+2)$-angle and a Serre functor in $(d+2)$-angulated categories.

\begin{theorem}\label{thm1}\emph{\cite[Theorem 4.5]{z1}}
Let $\C$ be a $(d+2)$-angulated category. Then $\C$ has Auslander-Reiten $(d+2)$-angles if and only if $\C$ has a Serre functor.
\end{theorem}

\section{Grothendieck groups and Auslander-Reiten $(d+2)$-angles}
In this section, let $k$ be an algebraically closed field. We always assume that $\C$ is a $k$-linear Hom-finite Krull-Schmidt $(d+2)$-angulated category.
We denote by $\ind(\C)$ the set of isomorphism classes of indecomposable objects in $\C$.
For any $X\in\ind(\C)$, we denote by $\supp \Hom_{\C}(X,-)$ the subcategory of $\C$ generated by
objects $Y$ in $\ind(\C)$ with $\Hom_{\C}(X,Y)\neq0$. Similarly, $\supp\Hom_{\C}(-,X)$ denotes
the subcategory generated by objects $Y$ in $\ind(\C)$ with $\Hom_{\C}(Y,X)\neq 0$. If
$\supp\Hom_{\C}(X,-)$ ($\supp\Hom_{\C}(-,X)$, respectively) contains only finitely many
indecomposables, we say that $|\supp\Hom_{\C}(X,-)|<\infty $ ($|\supp\Hom_{\C}(-,X)|<\infty$
respectively).

\begin{definition}\cite[Definition 3.1]{z1}
A $(d+2)$-angulated category $\C$ is called \emph{locally finite}
if $|\supp\Hom_{\C}(X,-)|<\infty$ and $|\supp\Hom_{\C}(-,X)|<\infty$, for any object $X\in\ind(\C)$.
\end{definition}

\begin{theorem}\emph{\cite[Theorem 3.8]{z1}}
Let $\C$ be a locally finite $(d+2)$-angulated category. Then
$\C$ has Auslander-Reiten $(d+2)$-angles.
\end{theorem}

The following result is essentially already in \cite[Lemma 2.1]{l}.

\begin{lemma}\label{m2}
Let $(\C,\Sigma^d,\Theta)$ be a $(d+2)$-angulated category, and
$$A_{\bullet}:~A_0\xrightarrow{\alpha_0}A_1\xrightarrow{\alpha_1}\cdots
\xrightarrow{\alpha_{d-3}}A_{d-2}\xrightarrow{\alpha_{d-2}}A_{d-1}\xrightarrow{\binom{u}{v}}L\oplus M\xrightarrow{(c,~d)}N\xrightarrow{~h~}\Sigma^d A_0$$
a $(d+2)$-angle in $\C$.
If $u=0$, then it is isomorphic to the following $(d+2)$-angle:
$$A_0\xrightarrow{\alpha_0}A_1\xrightarrow{\alpha_1}\cdots
\xrightarrow{\alpha_{d-3}}A_{d-2}\xrightarrow{\alpha_{d-2}}A_{d-1}\xrightarrow{\binom{0}{v}}L\oplus M \xrightarrow{\left(\begin{smallmatrix}
1&0\\
0&t
\end{smallmatrix}\right)}L\oplus Q\xrightarrow{~~}\Sigma^d A_0.
$$
Moreover, the $(d+2)$-angle $A_{\bullet}$ is a direct sum of the following two $(d+2)$-angles
$$B_{\bullet}:~A_0\xrightarrow{\alpha_0}A_1\xrightarrow{\alpha_1}\cdots
\xrightarrow{\alpha_{d-3}}A_{d-2}\xrightarrow{\alpha_{d-2}}A_{d-1}\xrightarrow{~v~}M\xrightarrow{~t~}Q\xrightarrow{~~}\Sigma^d A_0,$$
$$M_{\bullet}:~0\xrightarrow{~}0\xrightarrow{~}\cdots
\xrightarrow{~}0\xrightarrow{~}0\xrightarrow{~}L\xrightarrow{1_L}L\xrightarrow{~~}0.$$
\end{lemma}

\proof Since $u=0$, by (N3), we have the following commutative diagram
$$\xymatrix{
A_0\ar[r]^{\alpha_1}\ar@{-->}[d]& A_1 \ar[r]^{\alpha_2}\ar@{-->}[d] & A_2 \ar[r]^{\alpha_3}\ar@{-->}[d]& \cdots \ar[r]&A_{d-1}\ar[r]^{\binom{u}{v}\quad}\ar[d]&L\oplus M\ar[r]^{\quad(c,~d)}\ar[d]^{(1,~0)}&N\ar[r]\ar@{-->}[d]^{c'} & \Sigma^d A_0 \ar@{-->}[d]\\
0\ar[r]&0 \ar[r] & 0 \ar[r]& \cdots \ar[r]&0\ar[r] & L \ar[r]^{1}&L\ar[r]& 0}$$
It follows that $c'(c,d)=(1,0)$ and then $c'c=1$. This shows that $c$ is a section.

Now we can assume that $A_{\bullet}$ is the following form
$$A_0\xrightarrow{\alpha_0}A_1\xrightarrow{\alpha_1}\cdots
\xrightarrow{\alpha_{d-3}}A_{d-2}\xrightarrow{\alpha_{d-2}}A_{d-1}\xrightarrow{\binom{0}{v}}L\oplus M\xrightarrow{\left(\begin{smallmatrix}
s&w\\
0&t
\end{smallmatrix}\right)}u(M)\oplus Q\xrightarrow{~(x,~y)~}\Sigma^d A_0$$
where $s$ is an isomorphism.
Then $(x,y)\left(\begin{smallmatrix}
s&w\\
0&t
\end{smallmatrix}\right)=0$ implies $x=0$ since $s$ is an isomorphism, and
$\left(\begin{smallmatrix}
s&w\\
0&t
\end{smallmatrix}\right)\left(\begin{smallmatrix}
0\\v
\end{smallmatrix}\right)=0$ implies $wv=0$.
Thus $(0,w)\binom{0}{v}=wv=0$. So there exists a morphism
$(a,b)\colon u(M)\oplus Q\to u(M)$ such that $(0,w)=(a,b)\left(\begin{smallmatrix}
s&w\\
0&t
\end{smallmatrix}\right).$ In particular, $as=0$ and $aw+bt=w$.
Since $s$ is an isomorphism, we have $a=0$ implies $w=bt$.
Hence we have a commutative diagram
$$\xymatrix@R=1.5cm{A_{d-1}\ar[r]^{\binom{0}{v}\quad}\ar@{=}[d]&L\oplus M\ar@{=}[d]
\ar[r]^{\left(\begin{smallmatrix}
1&0\\
0&t
\end{smallmatrix}\right)}\ar@{=}[d]&L\oplus Q\ar[d]^{\left(\begin{smallmatrix}
s&b\\
0&1
\end{smallmatrix}\right)}\ar[r]^{(0,y)}&\Sigma^dA_0\ar@{=}[d]\\
A_{d-1}\ar[r]^{\binom{0}{v}\quad}&L\oplus M\ar[r]^{\left(\begin{smallmatrix}
s&w\\
0&t
\end{smallmatrix}\right)}&L\oplus Q\ar[r]^{(0,y)}&\Sigma^dA_0}$$
implies that $A_{\bullet}\simeq B_{\bullet}\oplus M_{\bullet}$ since $\left(\begin{smallmatrix}
s&b\\
0&1
\end{smallmatrix}\right)$ is an isomorphism.
Since $\Theta$ is closed under direct summands and $A_{\bullet}\in\Theta$
we have $B_{\bullet}\in\Theta$, that is, it is a $(d+2)$-angle.  \qed

\begin{lemma}\label{m3}
Let $\C$ be a locally finite $(d+2)$-angulated category. Then for any object $X$
in $\C$, there exists a natural number $n$ {\rm(}$m$, respectively{\rm)} such that
${\rm rad}_{\C}^n(-,X)=0$
{\rm(}${\rm rad}_{\C}^m (X,-)=0$, respectively{\rm)}.
\end{lemma}

\proof This result was proved in \cite[Lemma 1.2]{xz} for the case that $\C$ is a triangulated category. But their proof
can be applied to the
context of a $(d+2)$-angulated category without any change.  \qed
\medskip

Suppose that $\C$ is an essentially small $(d+2)$-angulated
category, and let $F(\C)$ be the free abelian group on the set of isomorphism classes $\langle A\rangle$ of objects $A$ in $\C$. Given a $(d+2)$-angle
$$A_{\bullet}:~A_0\xrightarrow{~}A_1\xrightarrow{~}A_2\xrightarrow{~}\cdots\xrightarrow{~}A_d\xrightarrow{~}A_{d+1}\xrightarrow{~}\Sigma^d A_0$$
in $\C$, we denote the corresponding Euler relation in $F(\C)$ by $\chi(A_{\bullet})$, that is,
$$\chi(A_{\bullet}):=[A_0]-[A_1]+[A_2]+\cdots+(-1)^{d+1}[A_{d+1}].$$

\begin{definition}\cite[Definition 2.1]{bt2} and \cite[Definition 2.2]{fe1}
Let $\C$ be an essentially small $(d+2)$-angulated category, and $F(\C)$
  the free abelian group on the set of isomorphism classes $[A]$ of objects $A$ in $\C$.  Morever, let $R(\C)$ be the
  subgroup of $F(\C)$ generated by the following sets of elements
 $$\{\chi(A_\bullet) \mid A_\bullet ~\text{is a $(d+2)$-angle in } \C \}$$
  in $\C$.  The \emph{Grothendieck group} $K_0(\C)$ of $\C$ is the
  quotient group $F(\C)/R(\C)$.  Given an object $A \in \C$, the
  residue class $\langle A \rangle + R(\C)$ in $K_0(\C)$ is denoted by
  $[A]$.
\end{definition}

\begin{definition}
Let $\C$ be an essentially small $(d+2)$-angulated category, and $F(\C)$
  the free abelian group on the set of isomorphism classes $\langle A
  \rangle$ of objects $A$ in $\C$.  Morever, let $R'(\C)$ be the
  subgroup of $F(\C)$ generated by the following sets of elements
   $$\{\chi(A_\bullet) \mid A_\bullet ~\text{is a \textbf{split} $(d+2)$-angle in } \C \}$$
  in $\C$.  The \emph{{\bf\emph{split}} Grothendieck group} $K_0(\C,0)$ of $\C$ is the
  quotient group $F(\C)/R'(\C)$.  Given an object $A \in \C$,  the
  residue class $\langle A \rangle + R'(\C)$ in $K_0(\C,0)$ is denoted by
  $[A]$. For simplicity, sometimes we also denote by $[A_{\bullet}]$ the element in $K_0(\C,0)$.
\end{definition}

Note that the definition of the Grothendieck group is the reason why we are only considering essentially small categories: the
collection of isomorphism classes in the category must form a set. In addition, there exists a canonical epimorphism $\phi\colon K_0(\C,0)\to K_0(\C)$.

\begin{remark}\cite[Proposition 2.2]{bt2}
Let $\C$ be an essentially small $(d+2)$-angulated category, and $K_0(\C,0)$ its split Grothendieck group.

(1) The element $[0]$ is the zero element in $K_0(\C,0)$.

(2) If $A$ and $B$ are objects in $\C$, then $[A\oplus B]=[A]+[B]$ in $K_{0}(\C,0)$.
\end{remark}

Now we are ready to state and prove our first main result.

\begin{theorem}\label{main1}
Let $\C$ be a locally finite $(d+2)$-angulated category. Then
$\Ker\phi$ is generated  by the elements $[A_{\bullet}]$ in $K_0(\C,0)$, where
$$A_{\bullet}:~A_0\xrightarrow{\alpha_0}A_1\xrightarrow{\alpha_1}A_2\xrightarrow{\alpha_2}\cdots\xrightarrow{\alpha_{d-1}}A_d\xrightarrow{\alpha_d}A_{d+1}\xrightarrow{\alpha_{d+1}}\Sigma^d A_0.$$
runs
through all Auslander-Reiten $(d+2)$-angles in $\C$.
\end{theorem}

\proof
 Let $A_{\bullet}:~A_0\xrightarrow{\alpha_0}A_1\xrightarrow{\alpha_1}A_2\xrightarrow{\alpha_2}\cdots\xrightarrow{\alpha_{d-1}}A_d\xrightarrow{\alpha_d}A_{d+1}\xrightarrow{\alpha_{d+1}}\Sigma^d A_0 $ be an arbitrary $(d+2)$-angle with
$\alpha_{d+1}\neq 0$ and $A_{d+1}\in\ind(\C)$. It suffices to prove that the element $[A_{\bullet}]$ in $K_0(\C,0)$ can be
written as a sum of the elements in $K_0(\C,0)$ corresponding to some Auslander-Reiten $(d+2)$-angles.

Suppose $\alpha_{d+1}\in{\rm rad}_{\C}^n(A_{d+1},\Sigma^d A_0)$, and
$$B_{\bullet}:~B_0\xrightarrow{\beta_0}B_1\xrightarrow{\beta_1}B_2\xrightarrow{\beta_2}\cdots\xrightarrow{\beta_{d-1}}B_d\xrightarrow{\beta_d}A_{d+1}\xrightarrow{\beta_{d+1}}\Sigma^d B_0 $$
is an Auslander-Reiten $(d+2)$-angle
ending at $A_{d+1}$. Since $\alpha_{d+1}\neq 0$, by Lemma \ref{y1}, we know that $\alpha_d$ is not a retraction.
Then there exists a morphism $\varphi_d\colon A_d\to B_d$ such that $\alpha_d=\beta_d\varphi_d$.
By Lemma \ref{m4}, there exists a morphism
$$\xymatrix{
A_0\ar[r]^{\alpha_0}\ar@{-->}[d]^{\varphi_{0}}& A_1 \ar[r]^{\alpha_1}\ar@{-->}[d]^{\varphi_{1}} & A_2 \ar[r]^{\alpha_2}\ar@{-->}[d]^{\varphi_{2}}& \cdots \ar[r]&A_{d-1}\ar[r]^{\alpha_{d-1}}\ar@{-->}[d]^{\varphi_{d-1}}&A_{d}\ar[r]^{\alpha_d}\ar[d]^{\varphi_d}&A_{d+1}\ar[r]^{\alpha_{d+1}}\ar@{=}[d] & \Sigma^d A_0 \ar@{-->}[d]^{\Sigma^d\varphi_0}\\
B_0\ar[r]^{\beta_0}&B_1 \ar[r]^{\beta_1} & B_2 \ar[r]^{\beta_2}& \cdots \ar[r]&B_{d-1}\ar[r]^{\beta_{d-1}} & B_{d} \ar[r]^{\beta_d}&A_{d+1}\ar[r]^{\beta_{d+1}}& \Sigma^dB_0}$$
of $(d+2)$-angles, moreover,
$$C_{\bullet}:~A_0\xrightarrow{\left(
                    \begin{smallmatrix}
                      -\alpha_1 \\
                      \varphi_0 \\
                    \end{smallmatrix}
                  \right)} A_1\oplus B_0\xrightarrow{\left(
                             \begin{smallmatrix}
                               \alpha_1 & 0 \\
                               \varphi_1 & \beta_0
                             \end{smallmatrix}
                           \right)}
 A_2\oplus B_1\xrightarrow{\left(
                                       \begin{smallmatrix}
                                         -\alpha_2 & 0  \\
                                         \varphi_2 &\beta_1
                                       \end{smallmatrix}
                                     \right)}\cdots\xrightarrow{\left(
                                       \begin{smallmatrix}
                                         (-1)^{d}\alpha_{d-1}& 0  \\
                                         \varphi_{d-1} &\beta_{d-2}
                                       \end{smallmatrix}
                                     \right)}A_d\oplus B_{d-1}$$
                                     $\xrightarrow{\left(
                                       \begin{smallmatrix}
                                         (-1)^{d+1}\varphi_d,&\beta_{d-1}
                                       \end{smallmatrix}
                                     \right)}B_d\xrightarrow{~\alpha_{d+1}\beta_d~}\Sigma^dA_0$
 is a $(d+2)$-angle in $\C$ with $\alpha_{d+1}\beta_d\in{\rm rad}_{\C}^{n+1}(B_d,\Sigma^d A_0)$, and
 that $[A_{\bullet}]=[B_{\bullet}]+[C_{\bullet}]$.

Put $C_i:=A_{i}\oplus B_{i-1},~i=1,2,\cdots,d$. At this time $C_{\bullet}$ is of the form
$$A_0\xrightarrow{~~}C_1\xrightarrow{~~}C_2\xrightarrow{~~}\cdots\xrightarrow{~~}C_d\xrightarrow{~~}B_{d}\xrightarrow{\alpha_{d+1}\beta_d}\Sigma^d A_0.$$
We decompose $B_d$
as a direct sum of indecomposable objects: $B_d= M_1\oplus M_2\oplus\cdots\oplus M_k$. Without
loss of generality, we can assume $k=2$. Then the $(d+2)$-angle $C_{\bullet}$ can be written as
$$A_0\xrightarrow{~~}C_1\xrightarrow{~~}C_2\xrightarrow{~~}\cdots\xrightarrow{~~}C_{d-1}\xrightarrow{~~}C_d\xrightarrow{~\binom{c}{e}~}M_1\oplus M_2\xrightarrow{(x,~y)}\Sigma^d A_0$$
with $(x,y)=\alpha_{d+1}\beta_d\in{\rm rad}_{\C}^{n+1}(B_d,\Sigma^d A_0)$.
 If $c$ is a retraction, we can assume that
$$C_{\bullet}:~A_0\xrightarrow{~~}C_1\xrightarrow{~~}C_2\xrightarrow{~~}\cdots\xrightarrow{~~}C_{d-1}\xrightarrow{~\binom{u}{v}~}M_1\oplus M'\xrightarrow{~\left(
                                       \begin{smallmatrix}
                                        1&0\\
                                        \ast&\ast
                                       \end{smallmatrix}
                                     \right)~}M_1\oplus M_2\xrightarrow{(x,~y)}\Sigma^d A_0$$
It follows that $\left(\begin{smallmatrix}
                                        1&0\\
                                        \ast&\ast
                                       \end{smallmatrix}
                                     \right)\left(
                                       \begin{smallmatrix}
                                        u\\v
                                       \end{smallmatrix}
                                     \right)=0$ and then $u=0$.
By Lemma \ref{m2}, the $(d+2)$-angle $C_{\bullet}$ is isomorphic to the
$(d+2)$-angle
$$D_{\bullet}:~A_0\xrightarrow{~~}C_1\xrightarrow{~~}C_2\xrightarrow{~~}\cdots\xrightarrow{~~}C_{d-1}\xrightarrow{~\binom{0}{v}~}M_1\oplus M'\xrightarrow{~\left(
                                       \begin{smallmatrix}
                                        1&0\\
                                        0&t
                                       \end{smallmatrix}
                                     \right)~}M_1\oplus Q\xrightarrow{(0,~z)}\Sigma^d A_0,$$
and $D_{\bullet}$ is a direct sum of the following $(d+2)$-angles:
$$N_{\bullet}:~A_0\xrightarrow{~}C_1\xrightarrow{~}\cdots
\xrightarrow{~}C_{d-2}\xrightarrow{~}C_{d-1}\xrightarrow{~v~}M'\xrightarrow{~t~}Q\xrightarrow{~z~}\Sigma^d A_0,$$
$$0\xrightarrow{~}0\xrightarrow{~}\cdots
\xrightarrow{~}0\xrightarrow{~}0\xrightarrow{~}M_1\xrightarrow{1_{M_1}}M_1\xrightarrow{~~}0.$$
Hence $z\in{\rm rad}_{\C}^{n+1}(Q,\Sigma^dA_0)$ since $(0,z)\in{\rm rad}_{\C}^{n+1}(Q,\Sigma^dA_0)$ and $[C_{\bullet}]=[N_{\bullet}]$.
Thus we can continue this process with $[N]_{\bullet}$ instead of $[C_{\bullet}]$. By Lemma \ref{m3}, we know that
this process must stop at a finite steps,
that is, up to some finite step, we can get a splitting $(d+2)$-angle. Then there are finitely
many Auslander-Reiten $(d+2)$-angles $Q_{\bullet}^{1},Q_{\bullet}^{2},\cdots,Q_{\bullet}^{q}$, such that
$[C_{\bullet}] =[N_{\bullet}]=[Q_{\bullet}^{1}]+[Q_{\bullet}^{2}]+\cdots+[Q_{\bullet}^{q}]$.
Hence we have proved the assertion in this case.
\medskip

If $e$ is a retraction, the proof of this result is similar to the case of $c$ is a retraction.

\medskip
Now we return to the $(d+2)$-angle $C_{\bullet}$ and assume that $c$ and $d$ are not
retraction in the following.  Assume that the Auslander-Reiten $(d+2)$-angle endings at $M_1, M_2$ are
respectively $U_{\bullet},V_{\bullet}$:
$$U_{\bullet}:~U_0\xrightarrow{~u_0~}U_1\xrightarrow{~u_1~}U_2\xrightarrow{~u_2~}\cdots\xrightarrow{~u_{d-1}~}U_d\xrightarrow{~u_d~}M_1
\xrightarrow{u_{d+1}}\Sigma^d U_0,$$
$$V_{\bullet}:~V_0\xrightarrow{~v_1~}V_1\xrightarrow{~v_2~}V_2\xrightarrow{~v_3~}\cdots\xrightarrow{~v_{d-1}~}V_d\xrightarrow{~v_d~}M_2
\xrightarrow{v_{d+1}}\Sigma^d V_0.$$
We form the direct sum of them: $U_{\bullet}\oplus V_{\bullet}$
$$U_0\oplus V_0\xrightarrow{~\delta_0~}U_1\oplus V_1\xrightarrow{~\delta_1~}U_2\oplus V_2\xrightarrow{\delta_2}\cdots\xrightarrow{\delta_{d-1}}U_d\oplus V_d\xrightarrow{\delta_d}M_1\oplus M_2
\xrightarrow{\delta_{d+1}}\Sigma^d U_0\oplus\Sigma^dV_0,$$
where $\delta_i=\left(\begin{smallmatrix}
 u_i&0\\
 0&v_i
 \end{smallmatrix}
 \right),~i=1,2,\cdots,d+1$.
Since $c$ is not retraction and $U_{\bullet}$ is an Auslander-Reiten $(d+2)$-angle, there exists a morphism $w_d\colon C_d\to U_d$
such that $c=u_dw_d$.
Since $e$ is not retraction and $V_{\bullet}$ is an Auslander-Reiten $(d+2)$-angle, there exists a morphism $w_d\colon C_d\to U_d$
such that $e=v_dw'_d$. It follows that $\left(\begin{smallmatrix}
 u_d&0\\
 0&v_d
 \end{smallmatrix}
 \right)\left(\begin{smallmatrix}
 w_d\\
 w'_d
 \end{smallmatrix}
 \right)=\left(\begin{smallmatrix}
 c\\
e
 \end{smallmatrix}
 \right)$.  By Lemma \ref{m4}, we have the following commutative diagram
$$\xymatrix{
A_0\ar[r]\ar@{-->}[d]& C_1 \ar[r]\ar@{-->}[d]& \cdots \ar[r]&C_{d-1}\ar[r]\ar@{-->}[d]&C_{d}\ar[r]^{\binom{c}{e}\qquad}\ar[d]^{\binom{u_d}{u'_d}}&M_1\oplus M_2\ar[r]^{~(x,~y)}\ar@{=}[d] & \Sigma^d A_0 \ar@{-->}[d]\\
U_0\oplus V_0\ar[r]^{\delta_0}&U_1\oplus V_1 \ar[r]^{\quad\delta_1}&\cdots \ar[r]&U_{d-1}\oplus
V_{d-1}\ar[r]^{\;\;\delta_{d-1}} & U_{d}\oplus V_d \ar[r]^{\delta_d\;\;}&M_1\oplus M_2\ar[r]^{\delta_{d+1}\;\;\;}& \Sigma^dU_0\oplus\Sigma^dV_0}$$
of $(d+2)$-angles such that
$$L_{\bullet}:~A_0\xrightarrow{~}C_1\oplus U_0\oplus V_0\xrightarrow{~}\cdots\xrightarrow{~}C_d\oplus
U_{d-1}\oplus V_{d-1}\xrightarrow{~}U_d\oplus V_d\xrightarrow{(x,y)\delta_d}\Sigma^dA_0$$
is a $(d+2)$-angle.
Since $(x,y)\in{\rm rad}_{\C}^{n+1}(B_d,\Sigma^d A_0)$, we have $(x,y)\delta_d\in{\rm rad}_{C}^{n+2}(U_d\oplus V_d,\Sigma^d A_0)$.
Thus we can continue this process with $[L]_{\bullet}$ instead of $[C_{\bullet}]$. By Lemma \ref{m3}, we know that
this process must stop at a finite steps,
that is, up to some finite step, we can get a splitting $(d+2)$-angle. Then there are finitely
many Auslander-Reiten $(d+2)$-angles $P_{\bullet}^{1},P_{\bullet}^{2},\cdots,P_{\bullet}^{p}$, such that
$[C_{\bullet}] =[P_{\bullet}^{1}]+[P_{\bullet}^{2}]+\cdots+[P_{\bullet}^{p}]$.
Hence we have proved the assertion in this case.

This completes the proof.  \qed

\begin{remark}
As a special case of Theorem \ref{main1} when $d=1$, it is just the Theorem 2.1 of Xiao and Zhu in \cite{xz}.
\end{remark}

\begin{example}
Let $\C$ be a triangulated category with a $d$-cluster titling subcategory $\X$ which is closed under the $d$-th
power of the shift functor. By \cite[Theorem 1]{gko}, we know that $\X$ is a $(d+2)$-angulated category.
If $\X$ is locally finite, by Theorem \ref{main1}, then the relations of its Grothendieck group are generated by all Auslander-Reiten
$(d+2)$-angulated category, see also \cite[Remark 5.5]{fe1}.
\end{example}

Now we  will show that the converse of Theorem \ref{main1} is also true.  For convenience, we will use the notation $[A,B]:=\dim_{k}\Hom_{\C}(A,B)$.

\begin{lemma}\label{x1}
Let
$A_{\bullet}:~A_0\xrightarrow{\alpha_0}A_1\xrightarrow{\alpha_1}A_2\xrightarrow{\alpha_2}\cdots\xrightarrow{\alpha_{d-1}}A_d\xrightarrow{\alpha_d}A_{d+1}\xrightarrow{\alpha_{d+1}}\Sigma^d A_0$ be an Auslander-Reiten $(d+2)$-angle in $\C$ and $U$ an object in $\C$. Then the following statements hold:
\begin{itemize}
\item[\rm(1)] The morphism ${\rm Hom}_{\C}(U,\alpha_d)\colon{\rm Hom}_{\C}(U,A_d)\xrightarrow{~}{\rm Hom}_{\C}(U,A_{d+1})$ is an epimorphism if and only if
$A_{d+1}$ is not a direct summand in $U$.

\item[\rm (2)] If $U\in\ind(\C)$, then the morphism ${\rm Hom}_{\C}(U,\alpha_0)\colon{\rm Hom}_{\C}(U,A_0)\xrightarrow{~}{\rm Hom}_{\C}(U,A_{1})$ is a monomorphism if and only if
$U\not\simeq\Sigma^dA_{d+1}$.

\item[\rm (3)] If $U \in \ind(\C)$, we have $[U,A_0]-[U,A_1]+[U,A_2]+\cdots+(-1)^{d+1}[U,A_{d+1}] \neq 0$ if and only if $U \simeq A_{d+1}$ or $U \simeq \Sigma^{-d}A_{d+1}$.
\end{itemize}
\end{lemma}

\proof (1) We know that $A_{d+1}$ is a direct summand in $U$ if and only if there exists a retraction  $U\to A_{d+1}$. By the definition of an Auslander-Reiten $(d+2)$-angle, this is equivalent to
${\rm Hom}_{\C}(U,\alpha_d)\colon{\rm Hom}_{\C}(U,A_d)\xrightarrow{~}{\rm Hom}_{\C}(U,A_{d+1})$ not being epimorphism, which
proves (1).

(2) Applying the functor $\Hom_{\C}(U,-)$ to the $(d+2)$-angle $A_{\bullet}$, we have the following exact sequence:
$$\cdots\xrightarrow{}\Hom_{\C}(U,\Sigma^{-d}A_{d})\xrightarrow{\beta}\Hom_{\C}(U,\Sigma^{-d}A_{d+1})
\xrightarrow{~~}\Hom_{\C}(U,A_0)\xrightarrow{~\gamma~}\Hom_{\C}(U,A_1)\xrightarrow{}\cdots$$
where $\beta=\Hom_{\C}(U,\Sigma^{-d}\alpha_{d+1})$ and $\gamma=\Hom_{\C}(U,\alpha_0)$.
It follows that the morphism $\Hom_{\C}(U,\alpha_0)$ is monomorphism if and only if the morphism
$\Hom_{\C}(U,\Sigma^{-d}\alpha_{d+1})$ is epimorphism.
By applying (1) to the object $\Sigma^{-d}U$, we obtain that
$\Hom_{\C}(U,\Sigma^{-d}\alpha_{d+1})$ is an epimorphism if and only if $A_{d+1}$ is not a direct summand in $\Sigma^dU$.
Since $U\in\ind{\C}$, we have that $\Sigma^dU$ is also indecomposable, this is equivalent to
$U\not\simeq \Sigma^{-d}A_{d+1}$. Thus (2) holds.

(3) Applying the functor $\Hom_{\C}(U,-)$ to the $(d+2)$-angle $A_{\bullet}$, we have the following exact sequence:
$$\cdots\xrightarrow{}\Hom_{\C}(U,A_0)\xrightarrow{~\mu~}\Hom_{\C}(U,A_1)
\xrightarrow{}\cdots\xrightarrow{}\Hom_{\C}(U,A_{d-1})\xrightarrow{~\nu~}\Hom_{\C}(U,A_{d+1})\xrightarrow{}\cdots$$
where $u=\Hom_{\C}(U,\alpha_0)$ and $\nu=\Hom_{\C}(U,\alpha_{d+1})$.
We also get the following exact sequence:
$$0\xrightarrow{}K\xrightarrow{}\Hom_{\C}(U,A_0)\xrightarrow{~\mu~}\Hom_{\C}(U,A_1)
\xrightarrow{}\cdots\xrightarrow{}\Hom_{\C}(U,A_{d-1})\xrightarrow{~\nu~}\Hom_{\C}(U,A_{d+1})\xrightarrow{}M\xrightarrow{}0$$
where $K=\Ker\mu$ and $M=\Coker\nu$.
Splitting into short exact sequences and using our finiteness assumption, we get that the
equation
$$[U,A_0]-[U,A_1]+[U,A_2]+\cdots+(-1)^{d+1}[U,A_{d+1}]=\dim_{k}K+\dim_{k}M$$
where $d$ is odd and
$$[U,A_0]-[U,A_1]+[U,A_2]+\cdots+(-1)^{d+1}[U,A_{d+1}]=\dim_{k}K-\dim_{k}M$$
where $d$ is even.
Hence $[U,A_0]-[U,A_1]+[U,A_2]+\cdots+(-1)^{d+1}[U,A_{d+1}]\neq 0$ if and only if the righthand sides of the two equations are also non-zero. This means that either $K$ or $M$ (or both) must be non-zero.
The object $K$ is non-zero if and only if $\Hom_{\C}(U,\alpha_0)$ is not monomorphism.
By (2), we know that $\Hom_{\C}(U,\alpha_0)$ is not monomorphism if and only if $U\simeq \Sigma^dA_{d+1}$.
Similarly, $M$ is non-zero if and only if $\Hom_{\C}(U,\alpha_{d+1})$ is not epimorphism.
By (1), we know that $\Hom_{\C}(U,\alpha_{d+1})$ is not epimorphism if and only if
$A_{d+1}$ is  a direct summand in $U$. Since $U$ is indecomposable, we have
$U\simeq A_{d+1}$.  This completes the proof.\qed

\begin{lemma}\label{y8}
Assume that $b_1[C_1]+b_2[C_2]+\cdots+b_m[C_m]=0$ in $K_0(\C,0)$ for integers
$b_i$ and objects $C_1,C_2,\cdots,C_m$ in $\C$.
Then we have $b_1[X,C_1]+b_2[X,C_2]+\cdots+b_m[X,C_m]=0$ in $\mathbb{Z}$ for any
object $X$ in $\C$.
\end{lemma}

\proof  Suppose $b_1[C_1]+b_2[C_2]+\cdots+b_m[C_m]=0$ in $K_0(\C,0)$. By moving negative
terms to the right-hand side of the equality, without loss of generality, we can
assume $b_i\geq 0$ for each $i =1,2,\cdots,m$. Using the defining Euler relations for
$K_0(\C,0)$, we obtain
$$b_1[C_1]+b_2[C_2]+\cdots+b_m[C_m]=[b_1C_1\oplus b_2C_2\oplus\cdots\oplus b_mC_m]=0$$
where $b_iC_i$ is the coproduct of the object $C_i$ with itself $b_i$ times.
Thus the object $b_1C_1\oplus b_2C_2\oplus\cdots\oplus b_mC_m$ is zero in $\C$.
Applying $[X,-]$ and using additivity, we get that
the equation $b_1[X,C_1]+b_2[X,C_2]+\cdots+b_m[X,C_m]=0$.  \qed

\begin{lemma}\label{x4}
Suppose that
$\Ker\phi$ is generated by the elements $[A_{\bullet}]$ in $K_0(\C,0)$, where
$$A_{\bullet}:~A_0\xrightarrow{\alpha_0}A_1\xrightarrow{\alpha_1}A_2\xrightarrow{\alpha_2}\cdots\xrightarrow{\alpha_{d-1}}A_d\xrightarrow{\alpha_d}A_{d+1}\xrightarrow{\alpha_{d+1}}\Sigma^d A_0.$$
runs through all Auslander-Reiten $(d+2)$-angles in $\C$.
If for any $X\in\ind(\C)$, there exists an indecomposable object $U$ in $\C$ such that
${\rm Hom}_{\C}(U,X)\neq 0$, then $\C$ is locally finite.
\end{lemma}

\proof Assume that for any $X\in\ind(\C)$, there exists an indecomposable object $U$ in $\C$ such that $\Hom_{\C}(U,X)\neq 0$.
Note that the split $(d+2)$-angle
$$O_{\bullet}:~\Sigma^{-d}X\to 0\to 0\to\cdots\to 0\to X\xrightarrow{1_X}X,$$
belongs to $\Ker\phi$. Thus there are  Auslander-Reiten $(d+2)$-angles
$B^i_{\bullet}$ such that $O_{\bullet}=\sum\limits_{i=1}^{r}a_iB^i_{\bullet}$.
where $$B^{i}_{\bullet}:~B^i_0\xrightarrow{\beta^i_0}B^i_1\xrightarrow{\beta^i_1}B^i_2\xrightarrow{\beta^i_2}\cdots\xrightarrow{\beta^i_{d-1}}B^i_d\xrightarrow{\beta^i_d}B^i_{d+1}\xrightarrow{\beta^i_{d+1}}\Sigma^d B^i_0.$$
By Lemma \ref{y8}, we obtain that the equality
$$[U,O_{\bullet}]=\sum\limits_{i=1}^{r}a_i[U,B^i_{\bullet}].$$
Since ${\rm Hom}_{\C}(U,X)\neq0$, we have $[U,O_{\bullet}]=\dim_{k}{\rm Hom}_{\C}(U,X)+\dim_{k}{\rm Hom}_{\C}(U,\Sigma^{-d}X)\neq0$.
Hence there exists an integer $i\in\{1,2,\cdots,r\}$
such that $[U,B^i_{\bullet}]\neq 0$, that is,
$$[U,B^i_0]-[U,B^i_1]+[U,B^i_2]+\cdots+(-1)^{d+1}[U,B^i_{d+1}] \neq 0.$$
By Lemma \ref{x1}, this means that the indecomposable
$U$ is isomorphic to an object in the finite set $\{B^i_{d+1},\Sigma^dB^i_{d+1}\}^r_{i=1}$.
This shows that $|\supp\Hom_{\C}(-,X)|<\infty$.

Since $\C$ has Auslander-Reiten $(d+2)$-angles, by Theorem \ref{thm1}, we know that $\C$ has a Serre functor $\mathbb{S}$.
It follows that there exists an isomorphism ${\rm Hom}_{\C}(X,\mathbb{S}U)\simeq D{\rm Hom}_{\C}(U,X)\neq0$
implies that $|\supp\Hom_{\C}(X,-)|<\infty$.
Therefore $\C$ is locally finite.  \qed
\medskip

Our second main result is the following.

\begin{theorem}\label{main2}
 Suppose that
$\Ker\phi$ is generated by the elements $[A_{\bullet}]$ in $K_0(\C,0)$, where
$$A_{\bullet}:~A_0\xrightarrow{\alpha_0}A_1\xrightarrow{\alpha_1}A_2\xrightarrow{\alpha_2}\cdots\xrightarrow{\alpha_{d-1}}A_d\xrightarrow{\alpha_d}A_{d+1}\xrightarrow{\alpha_{d+1}}\Sigma^d A_0.$$
runs
through all Auslander-Reiten $(d+2)$-angles in $\C$.
Then  $\C$ is locally finite.
\end{theorem}

\proof By Lemma \ref{x4}, it suffices to show that for any $X\in\ind(\C)$, there exists an indecomposable object $U$ in $\C$ such that
${\rm Hom}_{\C}(U,X)\neq 0$.

Since $\C$ has Auslander-Reiten $(d+2)$-angles, by Theorem \ref{thm1}, we kow that $\C$ has a Serre functor $\mathbb{S}$.
Because $X$ is indecomposable, we have that $\mathbb{S}^{-1}X$ is also indecomposable.
Thus there exists an Auslander-Reiten $(d+2)$-angle:
$$A_0\xrightarrow{~}A_1\xrightarrow{~}A_2\xrightarrow{~}\cdots\xrightarrow{~}A_d\xrightarrow{~}\mathbb{S}^{-1}X\xrightarrow{~h~}\Sigma^d A_0.$$
It follows that $h\colon \mathbb{S}^{-1}X\to \Sigma^d A_0$ is non-zero and $A_0$ is indecomposable.
By the isomorphism $0\neq h\in{\rm Hom}_{\C}(\mathbb{S}^{-1}X,\Sigma^dA_0)\simeq D{\rm Hom}_{\C}(\Sigma^dA_0,X)$,
we obtain ${\rm Hom}_{\C}(\Sigma^dA_0,X)\neq 0$.
Since $A_0$ is indecomposable, we get that $\Sigma^dA_0$ is also indecomposable.  \qed

\begin{remark}
As a special case of Theorem \ref{main2} when $d=1$, we get that the converse of Theorem 2.1 of Xiao and Zhu in \cite{xz} is also true.
\end{remark}

\section*{Acknowledgement}
The author is grateful to Osamu Iyama to point out a shortcoming in a previous version of this
article, and give me helpful advice.  I also would like to thank Francesca Fedele for many useful comments.

\textbf{Panyue Zhou}\\
College of Mathematics, Hunan Institute of Science and Technology, 414006, Yueyang, Hunan, P. R. China.\\
and \\
D\'{e}partement de Math\'{e}matiques, Universit\'{e} de Sherbrooke, Sherbrooke,
Qu\'{e}bec J1K 2R1, Canada.\\
E-mail: \textsf{panyuezhou@163.com}

\end{document}